# THE DISTRIBUTION OF ORDER STATISTICS UNDER SAMPLING WITHOUT REPLACEMENT

B. O'Neill,[*] *Deloitte Australia*[**]

**Abstract**

This paper examines the distribution of order statistics taken from simple-random-sampling without replacement (SRSWOR) from a finite population with values $1, ..., N$. This distribution is a shifted version of the beta-binomial distribution, parameterised in a particular way. We derive the distribution and show how it relates to the distribution of order statistics under IID sampling from a uniform distribution over the unit interval. We examine properties of the distribution, including moments and asymptotic results. We also generalise the distribution to sampling without replacement of order statistics from an arbitrary finite population. We examine the properties of the order statistics for inference about an unknown population size (called the German tank problem) and we derive relevant estimation results based on observation of an arbitrary set of order statistics. We also introduce an algorithm that simulates sampling without replacement of order statistics from an arbitrary finite population without having to generate the entire sample.

SIMPLE-RANDOM-SAMPLING WITHOUT REPLACEMENT; ORDER STATISTICS; RANK-ORDER STATISTICS; FPOS DISTRIBUTION; BETA-BINOMIAL DISTRIBUTION; GERMAN TANK PROBLEM.

## 1. The finite-population-order-statistic (FPOS) distribution

Sampling problems occur widely in statistical practice and simple-random-sampling without replacement (SRSWOR) is the core method of sampling from a population. It is useful to examine the behaviour of order statistics arising from SRSWOR from a finite population of the values $1, ..., N$. The order statistics from such a population can be used for quantile analysis or estimation of an unknown population size.

Some basic information on the distribution of order statistics is available in Wilks (1962), with some further properties given as exercises for the reader (pp. 243-245, 251-252). It is also mentioned in Arnold *et al* (1992, p. 54) and Evans *et al* (2006, pp. 20-21), without elaboration on the properties of the distribution. The topic has not received much more attention in the literature. Indeed, it may surprise some readers to learn that even detailed statistical texts on order statistics do not usually cover the topic of order statistics under SRSWOR, focussing instead on order statistics from IID sampling from a continuous distribution.

[*] E-mail address: ben.oneill@hotmail.com.
[**] Deloitte Australia, 8 Brindabella Circuit, Canberra ACT 2609.

In this paper we will examine the distribution of the order statistics and derive a number of the properties of this distribution. For completeness we will derive the initial properties of the distribution available in Wilks (1962), but we will also extend this analysis to look at a mixture result, asymptotic properties, sufficiency and ancilliarity properties, derivation of the joint distribution of a vector of order statistics, and extension to order statistics for an arbitrary finite population. We will also connect the distribution with well-known results pertaining to the distribution of order statistics from IID sampling from the continuous uniform distribution over the unit interval, giving a more holistic view of the connection between order statistics from SRSWOR from a finite population and order-statistics in IID samples.

Consider a finite population composed of elements $1, \ldots, N$ and suppose we take a sample of $1 \leq n \leq N$ values from this population using simple-random-sampling without replacement. We denote the sample values as $X_1, \ldots, X_n$ and the ordered sample as $X_{(1)} \leq \cdots \leq X_{(n)}$. We are interested in the distribution of the order statistic $X_{(k)}$ for a value $1 \leq k \leq n$. We will call this the *finite-population-order-statistic (FPOS) distribution* and we will denote its probability mass function by $\text{FPOS}(x|k, n, N) \equiv \mathbb{P}(X_{(k)} = x|n, N)$.

Following Wilks (1962, p. 243) and Arnold *et al* (1992, p. 54) we establish the mass function for the FPOS distribution using a combinatorial argument. There are $\binom{N}{n}$ possible samples that can be drawn from the population, with equal probability. To get the event $X_{(k)} = x$ we have one way of choosing this order statistic, and then $\binom{x-1}{k-1}$ ways of choosing the previous order statistics and $\binom{N-x}{n-k}$ ways of choosing the ensuing order statistics. Consequently, applying the multiplication principle of counting gives the probability mass function:

$$\text{FPOS}(x|k, n, N) = \frac{\binom{x-1}{k-1}\binom{N-x}{n-k}}{\binom{N}{n}} \qquad \text{for } x = k, \ldots, N - n + k.$$

The support of the distribution reflects the fact that there are $k - 1$ distinct sample values below the order statistic and $n - k$ distinct sample values above it. This distribution is a shifted version of the negative hypergeometric distribution or beta-binomial distribution, given by:

$$\text{FPOS}(x|k, n, N) = \text{NegHyper}(x - k|N, N - n, k)$$
$$= \text{BetaBin}(x - k|N - n, k, n - k + 1).$$

Both of these distributional connections can easily be established algebraically, but the framing in terms of the shifted negative hypergeometric distribution also has an intuitive explanation pertaining to the interpretation of that distribution arising from SRSWOR. Instead of fixing the elements in the population and then sampling $n$ values, suppose we do the process in reverse — i.e., we arbitrarily predesignate $n$ values in a population of size $N$ to be the "sample values" and then we select all the population elements in a random order via SRSWOR and label them $1, 2, \ldots, N$ as they are selected. It can be shown that this gives the same sampling distribution as the original method. Using this reverse method, suppose we pause our selection of elements from the population once there are $k$ of the predesignated "sample values" in our selection. Since the elements are labelled with increasing labels by their selection order, the last element selected is the order statistic $X_{(k)}$, which is equal to the total number of elements selected at the time we select the $k$th value that was in the predesignated sample. In the language of the negative hypergeometric distribution, we have a population of size $N$ with $N - n$ "successes" (those not in the sample) and we have sampled until there are $k$ "failures" at which point there are $X_{(k)} - k$ "successes". Using the standard logic of the negative hypergeometric distribution (see e.g., Johnson *et al* 2005, pp. 253-254), we therefore have:

$$\mathbb{P}(X_{(k)} = x | n, N) = \text{NegHyper}(x - k | N, N - n, k).$$

The connection to the shifted beta-binomial distribution is also useful in examining the FPOS distribution. Although we will establish it formally in Theorem 4 below, the connection to the beta-binomial immediately leads to the following well-known mixture form:

$$\text{FPOS}(x | k, n, N) = \int_0^1 \text{Bin}(x - k | N - n, u) \cdot \text{Beta}(u | k, n - k + 1) \, du.$$

Consequently, an order statistic from a finite population can be generated by the process:

$$X_{(k)} \sim k + \text{Bin}(N - n, U_{(k)}) \qquad U_{(k)} \sim \text{Beta}(k, n - k + 1).$$

Readers familiar with the theory of order statistics will recognise the distribution of $U_{(k)}$ as that of the $k$th order statistic generated from $n$ IID values from a uniform distribution over the unit interval. The present mixture result establishes an interesting connection to the distribution of an order statistic in the IID uniform case. We will explore this in greater detail later in the paper and discuss the relationship between order statistics in SRSWOR versus order statistics in IID sampling. Nevertheless, even with this preliminary result we can see that the case of order statistics from SRSWOR can be regarded as a kind of variation to this behaviour, where there is an extra step in the simulation of the relevant order statistics.

## 2. Properties of the FPOS distribution

Since the FPOS distribution is just a scaled and re-parameterised version of the beta-binomial distribution, the properties of the former are easy to determine from the known properties of the latter (see e.g., Tripathi, Gupta and Gurland 1994). In particular, the generating functions for the FPOS distribution are scaled versions of the ordinary hypergeometric function and the moments are closely related to the moments of the beta-binomial distribution. Although it is possible to appeal to known results for the beta-binomial distribution, in the present paper we derive distributional properties from scratch, to elucidate the mathematics of the distribution. Theorems 1-3 below show the rising factorial moments and the mean, variance, skewness and kurtosis of the distribution. Theorem 4 then formally establishes the mixture results we have previously mentioned.

**THEOREM 1 (Rising factorial moments):** Using the notation $x^{(r)} = x(x+1) \cdots (x+r-1)$ to denote the rising factorials, we have:

$$\mathbb{E}(X_{(k)}^{(r)}) = \frac{(N+1)^{(r)} \cdot k^{(r)}}{(n+1)^{(r)}}.$$

**THEOREM 2 (Mean and variance):** The mean and variance are:

$$\mathbb{E}(X_{(k)}) = \frac{N+1}{n+1} \cdot k \qquad \mathbb{V}(X_{(k)}) = \frac{(N+1)(N-n)}{(n+1)^2(n+2)} \cdot k(n-k+1).$$

**THEOREM 3 (Skewness and kurtosis):** The skewness is:

$$\mathbb{Skew}(X_{(k)}) = (n - 2k + 1)\left(1 + 2 \cdot \frac{N-n-1}{n+3}\right)\sqrt{\frac{n+2}{(N+1)(N-n)k(n-k+1)}}.$$

The kurtosis is:

$$\mathbb{Kurt}(X_{(k)}) = 3 + \frac{1}{N+1}\left[\frac{n(n+1)^3(n+2)}{(N-n)(n+3)(n+4)k(n-k+1)} - \frac{6(n+1)^2(n+2)}{(N-n)(n+3)(n+4)} + \frac{6(N+1)(n+1)^2(n+2)}{(n+3)(n+4)k(n-k+1)} - \frac{6(N+1)(5n+11)}{(n+3)(n+4)}\right].$$

**THEOREM 4 (Mixture characterisation):** The probability mass function can be written as:

$$\text{FPOS}(x|k,n,N) = \int_0^1 \text{Bin}(x-k|N-n,u) \cdot \text{Beta}(u|k,n-k+1)\,du\,.$$

The probability mass function and the moment results in Theorems 1-2 above are available in Wilks (1962, pp. 243-244) but we have derived them here for completeness. The skewness and kurtosis results in Theorem 3 are not available in other literature on the distribution at issue. The mixture characterisation in Theorem 4 is not given by Wilks, but it is now a well-known result pertaining to the beta-binomial distribution (indeed, it is often taken as the defining property of the distribution and is the reason for the name of that distribution). Once it is established that the FPOS distribution is a re-parameterised version of the shifted beta-binomial distribution, the mixture result naturally follows.

Some further intuition for the mixture representation in Theorem 4 is to frame the initial simple-random-sampling process via ranks of an underlying set of IID uniform random variables. We generate values $U_1, \ldots, U_N \sim$ IID $U(0,1)$ and use $\mathbf{U}_n = (U_1, \ldots, U_n)$ and $\mathbf{U}_N = (U_1, \ldots, U_N)$ to denote vectors with respective sizes equal to the sample and population size. We will then let $U_{(1:n)} < \cdots < U_{(n:n)}$ denote the order statistics in the sample and $U_{(1:N)} < \cdots < U_{(N:N)}$ denote the order statistics in the population. We generate the population $X_1, \ldots, X_N$ by taking the ranks:

$$X_i = \text{rank}(U_i, \mathbf{U}_N)\,.$$

(Taking $X_i = \text{rank}(U_i, \mathbf{U}_N)$ means that $U_i = U_{(X_i:N)}$ in the population order statistics.) The SRSWOR is obtained by taking the first $n$ values $X_1, \ldots, X_n$ from this population. The order statistic $X_{(k)}$ is the $k$th smallest value in the sample. This is equal to the number of values in the sample that are no greater than this value (which is always $k$), plus the number of values outside the sample where the underlying values used to generate the population are no greater than the value $U_{(k)}$. We can therefore write the order statistics in the sample as:

$$X_{(k)} = k + \sum_{i=n+1}^{N} \mathbb{I}(U_i \leq U_{(k)})\,.$$

Since the values in $\mathbf{U}_N$ are IID uniform values we can see that the summation term here is a sum of $N-n$ Bernoulli random variables with fixed probability $U_{(k)}$. Moreover, since $U_{(k)}$ is an order statistic from $n$ IID uniform random variables we have $U_{(k)} \sim \text{Beta}(k, n-k+1)$. This gives intuitive confirmation to the decomposition that comes from the mixture result.

Our remaining analysis looks at the exact and asymptotic behaviour of the order statistic and relates this to well-known results for order statistics in the context of IID continuous random variables. It is useful to express the order statistics in a form that scales to be comparable to the scale of the $U_{(k)}$ values. To do this, we define the scaled order statistics:

$$\widetilde{U}_{(k)} \equiv \frac{X_{(k)}}{N+1}.$$

It is then simple to establish that:

$$\mathbb{E}(\widetilde{U}_{(k)}) = \mathbb{E}(U_{(k)}) \qquad \mathbb{V}(\widetilde{U}_{(k)}) = \frac{N-n}{N+1} \cdot \mathbb{V}(U_{(k)}).$$

As we can see, this scaling yields discrete values $\widetilde{U}_{(k)}$ that are within the unit interval and have the same expectation as the order statistics for IID sampling from the uniform distribution over the unit interval. There is lower variance due to the forced gaps (of at least $1/(N+1)$) between the scaled order statistics when sampling without replacement. These forced gaps means that the order statistics are more evenly spread out than in the IID case, which results in "sandwiching" the order statistics within a narrower range, thereby reducing the variance. (Here we note that there is arguably a slight deficiency in the scaling which also contributes to this result. One could take the view that an appropriate scaling to obtain a discrete analogy to the uniform distribution would be to use the discrete points that are the midpoints of a partition of $N$ equally-spaced subsets of the unit interval. Our scaling is slightly more compacted than this, which also contributes to the lower variance.) We now give asymptotic forms for the moments and a formal limit result for the distribution when the parameters are large.

**THEOREM 5 (Asymptotic moments):** Consider any limit path of valid parameters for the FPOS distribution where $N \to \infty$, $n \to \infty$ and $k \to \infty$ in such a way that $n/N \to \lambda$ and $k/n \to \phi$ with fixed limits $0 < \lambda < 1$ and $0 < \phi < 1$. Under any limiting path of this kind, the asymptotic mean, variance, skewness and kurtosis are given respectively by:

$$\hat{\mu}_N = \phi \cdot N,$$

$$\hat{\sigma}_N^2 = \frac{1-\lambda}{\lambda} \cdot \phi(1-\phi) \cdot N,$$

$$\hat{\gamma}_N = \frac{2(\frac{1}{2}-\phi)(2-\lambda)}{\sqrt{\lambda(1-\lambda)\phi(1-\phi)}} \cdot \frac{2}{\sqrt{N}},$$

$$\hat{\kappa}_N = 3 + \frac{1}{N}\left[\frac{\lambda}{1-\lambda}\left(\frac{1}{\phi(1-\phi)} - 6\right) + \frac{6}{\lambda}\left(\frac{1}{\phi(1-\phi)} - 5\right)\right],$$

and the FPOS distribution is asymptotically unskewed and mesokurtic.

We will show that a standardised FPOS random variable converges in distribution to the normal under various stipulated limit conditions. To do this we will examine the moment generating function of the standardised version of $X_{(k)}$, which is given by:

$$m(t) \equiv \exp\left(t \cdot \frac{X_{(k)} - \mathbb{E}(X_{(k)})}{\mathbb{S}(X_{(k)})}\right).$$

**THEOREM 6 (Convergence in distribution):** Consider any limit path of valid parameters for the FPOS distribution where $N \to \infty$, $n \to \infty$ and $k \to \infty$ in such a way that $n/N \to \lambda$ and $k/n \to \phi$ with fixed limits $0 < \lambda < 1$ and $0 < \phi < 1$. Under any such limit path we have:

$$\lim_{N,n,k} m(t) = \exp(t^2/2).$$

**THEOREM 7 (Convergence in distribution):** Consider any limit path of valid parameters for the FPOS distribution where $N \to \infty$, $n \to \infty$ and $k \to \infty$ in such a way that $n/N \to 1$ and $k/N \to 0$ and $k(N-n)/N \to \infty$. Under any such limit path we have:

$$\lim_{N,n,k} m(t) = \exp(t^2/2).$$

**COROLLARY:** Under any limit paths in Theorems 6-7 the standardised FPOS random variable convergences in distribution to the standard normal distribution —i.e., we have:

$$\lim_{N,n,k} \mathbb{P}\left(\frac{X_{(k)} - \mathbb{E}(X_{(k)})}{\mathbb{S}(X_{(k)})} \leq z\right) = \Phi(z) \qquad \text{for all } z \in \mathbb{R}.$$

**COROLLARY:** Under any limit paths in Theorem 6 we have $\widetilde{U}_{(k)} \to \phi$ in probability.

The asymptotic result in Theorem 6 means that for large values of $N$ we can approximate the FPOS distribution by the normal distribution with reasonable accuracy. Since the latter is a continuous distribution there are many ways it can be formed to give a discrete approximation. One method is to use the normed-pointwise density approximation:

$$\widehat{\text{FPOS}}(x|k,n,N) = \frac{\text{N}(x|\mu_{k,n,N}, \sigma^2_{k,n,N})}{\sum_r \text{N}(r|\mu_{k,n,N}, \sigma^2_{k,n,N})} \qquad x = k, \ldots, N-n+k,$$

where $\mu_{k,n,N} = \mathbb{E}(X_{(k)})$ and $\sigma^2_{k,n,N} = \mathbb{V}(X_{(k)})$ are the mean and variance. It is quite useful to examine the accuracy of the normal approximation over a range of parameter values to get a sense of how large $N$ needs to be to get a good approximation. To do this we can compute the log-root-mean-squared error of the approximation, which is given by:

$$\mathrm{LRMSE}_N(k,n) \equiv \log\left(\sqrt{\frac{1}{N-n+1}\sum_{x=k}^{N-n+k}(\mathrm{FPOS}(x|k,n,N)-\widehat{\mathrm{FPOS}}(x|k,n,N))^2}\right).$$

In Figure 1 we show heatmaps of $\mathrm{LRMSE}_N$ for the population sizes $N = 100, 200, 500, 1000$ for all allowable values of $k$ and $n$. These heatmaps confirm that the accuracy of the normal approximation gets better (lower LRMSE) as $N$ increases. The region of high accuracy near the line $k = (n+1)/2$ occurs because the skewness of the FPOS distribution on this line is zero (and near to this line the skewness is near zero). Further away from this line the FPOS distribution is more skewed and so the normal approximation is less accurate.

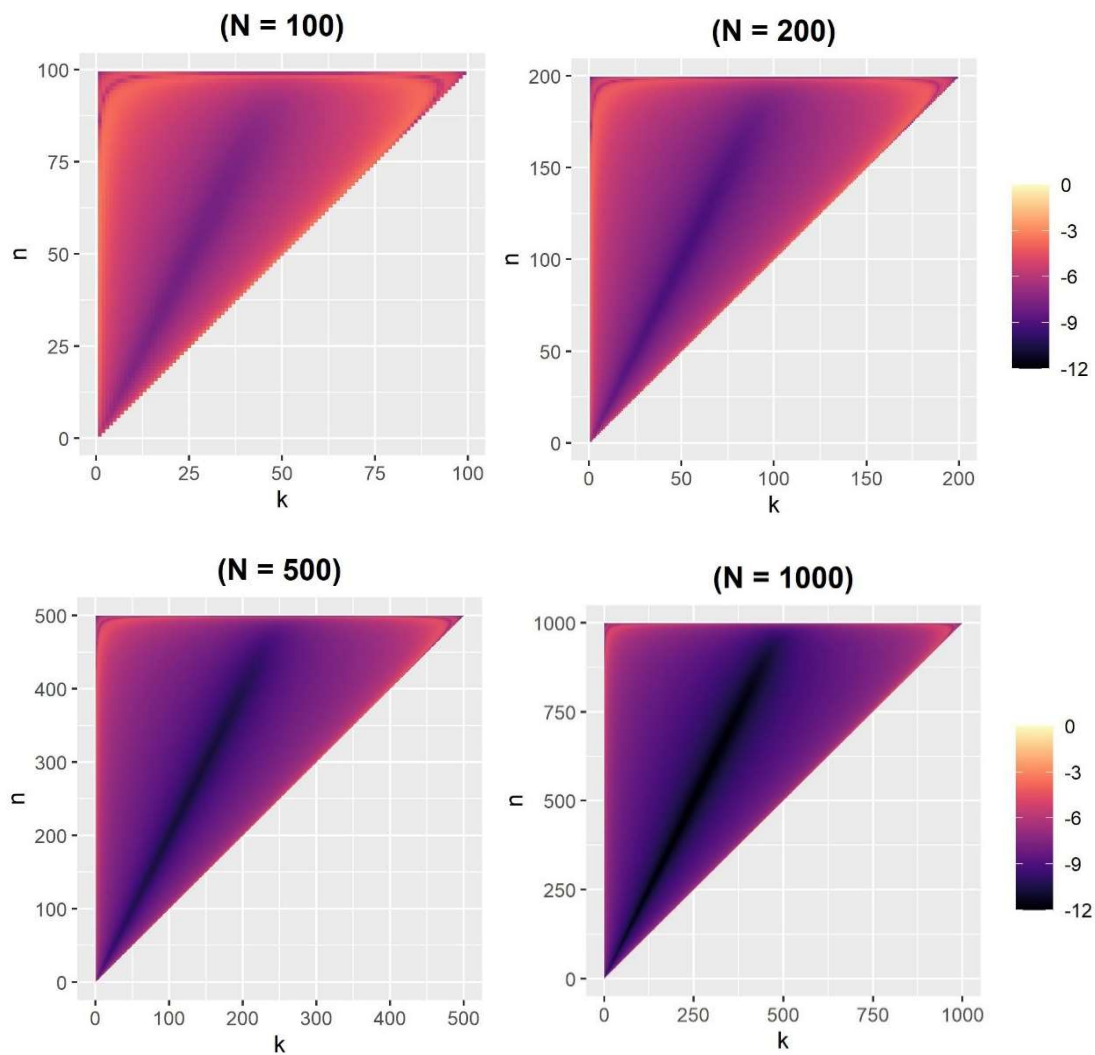

**FIGURE 1:** Heatmaps of the LRMSE of the normal approximation

## 3. Extension to the joint distribution of order statistics

Our previous analysis examined a single order statistic but we will now extend our analysis to look at the joint distribution of multiple order statistics. Suppose we now consider the vector of order statistics $\boldsymbol{X}_* \equiv (X_{(k_1)}, \ldots, X_{(k_r)})$ for some arbitrary ranks $1 \leq k_1 < \cdots < k_r \leq n$. To facilitate our analysis we also define the observed order statistics $\boldsymbol{x}_* \equiv (x_{(k_1)}, \ldots, x_{(k_r)})$ and the rank vector $\boldsymbol{k}_* \equiv (k_1, \ldots, k_r)$. As before, we can establish the probability mass function for the joint distribution of the order statistics using a combinatorial argument. This argument hinges on a simple probabilistic property of SRSWOR — if we condition on the event $X_{(k_i)} = x_i$ then the remaining sample values corresponding to higher order statistics are a SRSWOR from the elements $x_i + 1, \ldots, N$. Consequently, we have:

$$\mathbb{P}(X_{(k_{i+1})} = x | X_{(k_i)} = x_i) = \text{FPOS}(x - x_i | k_{i+1} - k_i, n - k_i, N - x_i).$$

Applying the multiplication rule of conditional probability then gives the mass function:

$$\text{FPOS}(\boldsymbol{x}_* | \boldsymbol{k}_*, n, N) = \frac{\binom{x_1-1}{k_1-1}\binom{x_2-x_1-1}{k_2-k_1-1}\binom{x_3-x_2-1}{k_3-k_2-1}\cdots\binom{x_r-x_{r-1}-1}{k_r-k_{r-1}-1}\binom{N-x_r}{n-k_r}}{\binom{N}{n}}.$$

This distribution is a shifted version of the Dirichlet-multinomial distribution, written in terms of differences pertaining to the order statistics and their ranks. In order to express the FPOS distribution in this form we define the difference vectors $\Delta_x \equiv \Delta_x(\boldsymbol{x}_*) \equiv (\Delta_{x,1}, \ldots, \Delta_{x,r}, \Delta_{x,r+1})$ and $\Delta_k \equiv \Delta_k(\boldsymbol{k}_*) \equiv (\Delta_{k,1}, \ldots, \Delta_{k,r}, \Delta_{k,r+1})$ by:

$$\Delta_{x,1} \equiv x_1 \qquad \Delta_{x,i+1} \equiv x_{i+1} - x_i \qquad \Delta_{x,r+1} \equiv N - x_r + 1,$$
$$\Delta_{k,1} \equiv k_1 \qquad \Delta_{k,i+1} \equiv k_{i+1} - k_i \qquad \Delta_{k,r+1} \equiv n - k_r + 1.$$

We can write the distribution as $\text{FPOS}(\boldsymbol{x}_* | \boldsymbol{k}_*, n, N) = \text{DirichletMu}(\Delta_x - \Delta_k | N - n, \Delta_k)$. This connection immediately leads to the following mixture form (Johnson *et al* 1997, pp. 200-207):

$$\text{FPOS}(\boldsymbol{x}_* | \boldsymbol{k}_*, n, N) = \int_\Theta \text{Mu}(\Delta_x - \Delta_k | N - n, \boldsymbol{s}) \cdot \text{Dirichlet}(\boldsymbol{s} | \Delta_k) \, d\boldsymbol{s},$$

where the space $\Theta$ is the probability space for the vector $\boldsymbol{s} = (s_1, \ldots, s_{r+1})$. Formally, this is the probability simplex $\Theta \equiv \{(s_1, \ldots, s_{r+1}) \in \mathbb{R}^{r+1} | \min_i s_i \geq 0, \sum_i s_i = 1\}$. Consequently, a vector of order statistics from a finite population can be generated by the process:

$$\boldsymbol{X}_* = \boldsymbol{k}_* + \Lambda \, \text{Mu}(N - n, \boldsymbol{S}) \qquad \boldsymbol{S} \sim \text{Dirichlet}(\Delta_k) \qquad \Lambda = \begin{bmatrix} 1 & 0 & \cdots & 0 & 0 \\ 1 & 1 & \cdots & 0 & 0 \\ \vdots & \vdots & \ddots & \vdots & \vdots \\ 1 & 1 & \cdots & 1 & 0 \end{bmatrix},$$

where $\boldsymbol{\Lambda}$ is an $r \times (r+1)$ matrix. As before, there is an interesting connection between the FPOS distribution and the distribution of order statistics in the IID uniform case. The vector $\boldsymbol{S}$ represents differences in the order statistics (at the ranks in $\boldsymbol{k}_*$) generated from $n$ IID samples from a uniform distribution over the unit interval, and the order statistics are obtained by the summation $\boldsymbol{U}_{(\boldsymbol{k})} = \boldsymbol{\Lambda} \boldsymbol{S}$. Although the Dirichlet-multinomial mixture is a known representation of the Dirichlet-multinomial distribution, we establish it here for completeness.

**THEOREM 8 (Mixture characterisation):** The mass function can be written as:

$$\text{FPOS}(\boldsymbol{x}_*|\boldsymbol{k}_*, n, N) = \int_\Theta \text{Mu}(\Delta_x - \Delta_k | N - n, \boldsymbol{s}) \cdot \text{Dirichlet}(\boldsymbol{s}|\Delta_k) \, d\boldsymbol{s},$$

where the space $\Theta$ is the probability space for the vector $\boldsymbol{s} = (s_1, \ldots, s_{r+1})$.

We do not examine the asymptotic behaviour of this distribution in the present paper, but it seems reasonable to conjecture that a standardised multivariate FPOS random variable would converge in distribution to the multivariate normal when all the parameters are large. This is a reasonable intuitive conjecture, since the Dirichlet random variable converges to a constant and the multinomial converges to the multivariate normal when the parameters are large.

We can use the joint distribution above to obtain useful properties of the order statistics relating to the population size $N$ when this is considered as an unknown parameter subject to inference. Intuition tells us that lower order statistics should not add any information about $N$ once any higher order statistic has been observed. Consequently, given multiple order statistics, the only one relevant for inference for $N$ is the largest. This is confirmed in Theorem 9, which shows that the highest observed order statistic is complete and sufficient for $N$ and the remaining order statistics are then conditionally ancillary for $N$.

**THEOREM 9 (Sufficient and ancillary statistics):** For any ranks $1 \leq k_1 < \cdots < k_r \leq n$ we consider the random vector $\boldsymbol{X}_* \equiv (X_{(k_1)}, \ldots, X_{(k_r)})$ and subvector $\boldsymbol{X}_{**} \equiv (X_{(k_1)}, \ldots, X_{(k_{r-1})})$. If we observe $\boldsymbol{X}_*$ as our sample then the following properties hold:
   (a) The statistic $X_{(k_r)}$ is complete and sufficient for $N$; and
   (b) The statistic $\boldsymbol{X}_{**}$ is ancillary for $N$ once we condition on $X_{(k_r)}$.

These two statistical properties establish that the highest observed order statistic is the only one relevant to inference about the population size. The underlying nature of SRSWOR means that once any order statistic is observed, any lower order statistics give no further information about the population size. The highest observed order statistic gives a lower bound to the population size but it is also a sufficient statistic for the population size. This fact is used in statistical estimation problems where $N$ is unknown.

**4. The German tank problem and extensions**

The inference problem of estimating $N$ is a famous statistical problem that is commonly known as the "German tank problem". The problem arose in WWII with American efforts to estimate the number of tanks in the German army based on the serial numbers of captured tanks (see Ruggles and Brodie 1947; Goodman 1952; Goodman 1954). In this historical military context, the German tanks were labelled by consecutive serial numbers $1, \ldots, N$ up to some unknown population size and the allies had access to a sample of $n$ captured tanks where the serial numbers had been inspected. Inference for $N$ was conducted by assuming that the sample was SRSWOR from the population. Ruggles and Brodie (1947) describe the available sample information in detail (see esp. pp. 78-79). The German production authorities had allocated manufacturers "blocks" of one hundred series numbers which were incompletely produced, with serial numbers occurring up to an unknown maximum. The estimation problem involved repeated estimation of the unknown number of tanks produced in each block, based on sample data from captured tanks with serial numbers in that block. Clark, Gonye and Miller (2021) consider an extension to the problem in which the population numbering does not start at one and so the goal is to estimate the population range based on observation of the smallest and largest order statistics.

One interesting aspect of the German tank problem is that —because it involves inference about an unknown population size $N$ based on a sample of size $n$— the size of the observable sample is restricted by the unknown population size. This is unlike many statistical inference problems where the sample size $n$ can be planned in advance. Indeed, if the sample size were planned in advance in this problem then it might not be possible to complete the planned sampling due to all the units in the population already having been observed prior to reaching the planned sample size. Instead, the German tank problem involves a case where new observations come in sporadically over time, without full control over the sample size. There is an ongoing effort to estimate the population size at any given time based on the sample size at that time. As a result, the assumption that sampling is SRSWOR should be considered carefully, having due regard to the fact that observation of sample items may occur in complicated circumstances.

The standard form of the German tank problem is an inference where we estimate $N$ using the highest order statistic $X_{(n)}$ (see also Gum *et al* 2005). We have $\mathbb{E}(X_{(n)}) = n(N+1)/(n+1)$ from Theorem 2 so we obtain an unbiased estimator for $N$ by scaling this order statistic as:

$$\widehat{N}_n \equiv \frac{n+1}{n} \cdot X_{(n)} - 1.$$

It is well-known that this estimator is the minimum variance unbiased estimator (MVUE) for the problem; this is a simple consequence of applying the Lehmann-Scheffé theorem (Casella and Berger 2001, p. 369) using Theorem 9. Using Theorem 2 and our asymptotic analysis, the variance of this estimator (and its asymptotic equivalent for large $n$ and $N$) are given by:

$$\mathbb{V}(\widehat{N}_n) = \frac{(N+1)(N-n)}{n(n+2)} \simeq \frac{1-\lambda}{\lambda^2}.$$

This estimator is the standard classical estimator used in inference problems of this kind, owing to its status as the MVUE. This estimator is widely applied when estimating the population size with sampling by SRSWOR. While the underlying assumption of SRSWOR was tenuous in the historical context of captured tanks in WWII, it is notable that this estimator performed well in the estimation of German tank numbers, famously outperforming alternative estimates from direct intelligence sources. (It is also interesting to note that the German army were victims of their own systematic approach to numbering of their tanks; the systematic nature of the serial numbers gave the allies useful intelligence which was exploited by statisticians. Countermeasures against this type of inference are now used in military applications, which reduce the information given by serial numbers on military resources that can be captured by an enemy.) Ruggles and Brodie (1947) observe that prior to the statistical estimation, "…Allied intelligence still suffered from the myth of German invincibility created by Nazi propagandists out of the successful *blitzkrieg* tactics in Poland and France, and it had grossly overestimated the enemy's position; the serial number technique revealed this fact and introduced realism in our picture of the strength of the German war machine" (pp. 80-81). The success of the statistical estimator under the suggests that it is somewhat robust to modest deviations from the assumed sampling method.

The German tank problem assumes that we observe an entire sample, so that we can obtain all order statistics including the highest order statistic (i.e., the maximum). Our previous results in Theorem 9 show that the highest observed order statistic is sufficient for $N$ and the remaining order statistics give no further information about this parameter. Of course, this assumes that the sample is taken by SRSWOR; the entire sample can be used to test this assumption if it is uncertain. Nevertheless, we can also consider the broader case where some order statistics are unobserved, generalising this analysis to consider how to estimate the population size if we observe an arbitrary order statistic $X_{(k)}$ without any of the higher order statistics. Here we assume that the "missingness" of the higher order statistics is not related to the values of those order statistics or any other order statistics (i.e., they are "missing completely at random"). If the unavailability of the higher order statistics is statistically related to the values of the order statistics then the analysis becomes far more complex and the results shown here do not hold.

In the case described, Theorem 9 ensures that the statistic $X_{(k)}$ is complete and sufficient for $N$. We can scale it to obtain the estimator:

$$\widehat{N}_k \equiv \frac{n+1}{k} \cdot X_{(k)} - 1,$$

which is an unbiased estimator with variance:

$$\mathbb{V}(\widehat{N}_k) = \frac{(N+1)(N-n)}{n+2} \cdot \frac{n-k+1}{k}.$$

This estimator remains the MVUE for the generalised version of the problem (again using the Lehmann-Scheffé theorem and Theorem 9). It can easily be shown that $\mathbb{V}(\widehat{N}_k)$ is a strictly decreasing function of $k$, so the variance is minimised when using the largest order statistic. This accords with our intuition that we can estimate the population size most accurately using the order statistic that is closest to the population size.

Some intuition for the estimator above is captured by supposing that the $n$ order statistics fall at equidistant points on the continuum $[0, N+1]$ (for this intuitive scenario we allow non-integer values for the order statistics) also with equidistance from the boundaries. This setup is shown in Figure 2 below. (It is notable that for the special case $n = N$ (i.e., a full census of the population) the order statistics would then fall on the values $1, \ldots, N$ which is the exact result.) Using this placement of points on the continuum, the $k$th order statistic is:

$$X_{(k)} = \frac{k}{n+1} \times (N+1).$$

Looking now at the formula for $\widehat{N}_k$ we can see that the estimator involves scaling up the order statistic by the multiplier $(n+1)/k$ to estimate the endpoint $N+1$ of this continuum. This scaling up involves inflating the order statistic in this setup up to the inferred place of the $n$th order statistic and then scaling it up one additional "unit" of distance to get to the inferred endpoint. In this placement of the order statistics we have $\widehat{N}_k = N$ so that the estimator for the population size is equal to the true population size. This illustrative approach and the corresponding intuition is similar to the spacing of quantiles in a QQ plot.

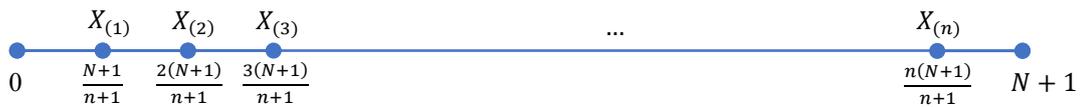

**FIGURE 2:** Order statistics falling at equidistant points on the continuum $[0, N+1]$

To examine the asymptotic behaviour of the population size estimator we can examine the ratio of the estimator to the true population size, which can be rewritten as:
$$\frac{\widehat{N}_k}{N} = \frac{n+1}{k} \cdot \frac{N+1}{N} \cdot \widetilde{U}_{(k)} - \frac{1}{N}.$$
Taking the limits $N \to \infty$, $n \to \infty$ and $k \to \infty$ and applying the corollary to Theorem 6 we have $\widehat{N}_k/N \to 1$ in probability. In this sense, the estimator $\widehat{N}_k$ converges towards the value $N$ (this is not the usual form of weak consistency, but it is useful nonetheless).

## 5. Relationship to the general distribution of order statistics in a finite population

In the above sections we have looked at order statistics for a population with values $1, \ldots, N$. More generally, in any finite population with values $z_1, \ldots, z_N$ (which may include non-distinct values) the rank-order statistics $r_1, \ldots, r_N$ are related to the order statistics $z_{(1)} \leq \cdots \leq z_{(N)}$ by the equation $z_k = z_{(r_k)}$, which lets us write the order statistics in terms of rank-order statistics and the ordered population values. It is quite simple to extend our above analysis to relate the FPOS distribution for the order statistics for a population with values $1, \ldots, N$ to the more general distribution for the order statistics of an arbitrary finite population.

For this analysis, we will consider $z_1, \ldots, z_N$ to be an arbitrary set of values that may include non-distinct values and we will relate the order statistics $z_{(1)} \leq \cdots \leq z_{(N)}$ back to the statistics $X_{(k)} = k$ using the mapping:

$$\omega(k) = z_{(k)}.$$

The mapping $\omega: \{1, \ldots, N\} \to \mathbb{R}$ is not generally injective so its inverse is a multivalued set function $\omega^{-1}$ defined as $\omega^{-1}(z) = \{x = 1, \ldots, N | \omega(x) = z\}$ for all $z \in \mathbb{R}$. Using this inverse function we can write the probability mass of the order statistic $Z_{(k)}$ as:

$$\mathbb{P}(Z_{(k)} = z) = \mathbb{P}(\omega(X_{(k)}) = z)$$
$$= \mathbb{P}(X_{(k)} \in \omega^{-1}(z))$$
$$= \sum_{x:\omega(x)=z} \text{FPOS}(z|k,n,N)$$
$$= \text{FPOS}(z|k,n,\omega),$$

where we define the following extension to the FPOS distribution:

$$\text{FPOS}(z|k,n,\omega) = \sum_{x:\omega(x)=z} \frac{\binom{x-1}{k-1}\binom{N-x}{n-k}}{\binom{N}{n}} \quad \text{for } z \in \{\omega(k), \ldots, \omega(N-n+k)\}.$$

This generalised version of the distribution allows for an arbitrary finite population where the values are not necessarily equally-spaced or distinct. The order statistics in this latter form can be generated by the random process:

$$Z_{(k)} = \omega(X_{(k)}) \qquad X_{(k)} \sim k + \text{Bin}(N-n, U_{(k)}) \qquad U_{(k)} \sim \text{Beta}(k, n-k+1).$$

Obviously, taking the mapping $\omega$ to be the identity function leads us back to the canonical form of the distribution with the population having elements $1, \ldots, N$. This present form is therefore a generalisation of the canonical form we have previously examined.

Following our analysis in Section 3, we can extend this characterisation of the order statistics to give a random process that generates vectors of order statistics. As before, we taken a set of arbitrary ranks $1 \leq k_1 < \cdots < k_r \leq n$ and we let $\boldsymbol{X}_* \equiv (X_{(k_1)}, \ldots, X_{(k_r)})$ denote the vector of rank-order statistics for these ranks. Similarly, we now let $\boldsymbol{Z}_* \equiv (Z_{(k_1)}, \ldots, Z_{(k_r)})$ denote the order statistics for these ranks and we say that $\boldsymbol{Z}_* = \boldsymbol{\omega}(\boldsymbol{X}_*)$ where $\boldsymbol{\omega}$ is the "vectorised" version of the mapping $\omega$. We can then generate the set of order statistics using the random process:

$$\boldsymbol{Z}_* = \boldsymbol{\omega}(\boldsymbol{X}_*) \qquad \boldsymbol{X}_* = \boldsymbol{k}_* + \boldsymbol{\Lambda}\operatorname{Mu}(N-n, \boldsymbol{S}),$$

$$\boldsymbol{S} \sim \operatorname{Di}(\boldsymbol{\Delta_k}) \qquad \boldsymbol{\Lambda} = \begin{bmatrix} 1 & 0 & \cdots & 0 & 0 \\ 1 & 1 & \cdots & 0 & 0 \\ \vdots & \vdots & \ddots & \vdots & \vdots \\ 1 & 1 & \cdots & 1 & 0 \end{bmatrix}.$$

For computing purposes, it is useful to note that the Dirichlet random variables in this process can be generated using IID gamma random variables (see e.g., Devroye 1986, pp. 593-599). This means we can "simplify" to a more primitive random process:

$$\boldsymbol{Z}_* = \boldsymbol{\omega}(\boldsymbol{X}_*) \qquad \boldsymbol{X}_* = \boldsymbol{k}_* + \boldsymbol{\Lambda}\operatorname{Mu}(N-n, \boldsymbol{S}),$$

$$S_i = \frac{R_i}{\sum_{\ell=1}^{r} R_\ell} \qquad R_i \sim \operatorname{Ga}(\Delta_{k,i}, 1),$$

$$\boldsymbol{\Lambda} = \begin{bmatrix} 1 & 0 & \cdots & 0 & 0 \\ 1 & 1 & \cdots & 0 & 0 \\ \vdots & \vdots & \ddots & \vdots & \vdots \\ 1 & 1 & \cdots & 1 & 0 \end{bmatrix}.$$

Since most statistical programming languages can generate gamma and multinomial pseudo-random values using built-in functions, it is relatively simple to program an algorithm that uses the above process to generate simulations of pseudo-random vectors of order statistics obtained from SRSWOR from a finite population with stipulated values. Below we give an algorithm in **R** to implement this random process to produce pseudo-random vectors of this kind. (This algorithm actually produces a matrix output where each row is one simulation.) This method avoids having to generate an entire sample to obtain the order statistics, so it constitutes an efficient method of simulation in cases where the number of elements in $\boldsymbol{k}_*$ is relatively small compared to the sample size $n$.

```
          Algorithm: Simulate vectors of order statistic
       using SRSWOR from an arbitrary finite population
```
| | |
|---|---|
| **Inputs:** | A vector of population values **population**<br>The sample size **size** (positive integer not larger than population)<br>A vector of rank-orders **ranks** (positive integers not larger than size)<br>Number of simulations **sims** to generate (positive integer) |
| **Output:** | A matrix with one row for each rank where each row is a simulation<br>of the vector of order statistics using the stipulated input values |

```
#Set sorted population and ranks
POP    <- sort(population)
```

```
RANKS  <- sort(ranks)
ORDER  <- order(ranks)

#Set parameters
N     <- length(POP)
r     <- length(RANKS)

#Set rank-difference vector
DELTA <- c(RANKS[1], diff(c(RANKS, size+1)))

#Generate OUT
OUT <- matrix(0, nrow = sims, ncol = r)
for (i in 1:sims) {
  R <- rgamma(r+1, shape = DELTA, scale = 1)
  S <- R/sum(R)
  V <- rmultinom(1, size = N-size, prob = S)
  X <- RANKS + cumsum(V[1:r])
  OUT[i, ] <- POP[X][ORDER] }

#Give output
OUT
```

In benchmark testing on this method we found that it performed significantly faster than the "standard method" that consists of generating the entire sample vector and extracting relevant order statistics through sorting. We tested both methods using the parameters $N = 40, n = 20$ and $\boldsymbol{k}_*$ composed of six elements, generating 1,000 simulations from each method 1,000 times. To measure computational time in a manner that is roughly invariant to machine-speed and background processes we included a third benchmark computation where we sort the numbers 1000:1 into ascending order (we call the time taken for this operation a "kilosort" and use it as a unit of measurement for the time taken for the two simulation methods). In our benchmarking test the algorithm used here took an average of 116.83 kilosorts whereas the standard algorithm took an average of 616.19 kilosorts. Our algorithm requires initial computation of the object $\Delta_{\boldsymbol{k}}$ to facilitate analysis, but after this it can simulate order statistics using a method that does not require generating the full sample. Consequently, it would be expected to be most efficient in cases where the number of ranks in the simulation is much smaller than the sample size; if the number of ranks for the order statistics is close to the sample size then it may be just as efficient to simulate the entire sample.

While the canonical form for the rank-order statistics has elegant and compact moment results, the general form gives the case-by-case result:

$$\mathbb{E}(f(Z_{(k)})) = \mathbb{E}(f \circ \omega(X_{(k)})) = \sum_{x=k}^{N-n+k} f(\omega(x)) \cdot \frac{\binom{x-1}{k-1}\binom{N-x}{n-k}}{\binom{N}{n}}.$$

This result does not simplify any further, owing to the generality of the mapping $\omega$. Similarly, there are no simple asymptotic results for $N \to \infty$ since such results would depend on the nature of the resulting sequence of population values that are used when taking the limit. If one is willing to stipulate properties of this sequence (e.g., its limiting distribution) then it is possible to obtain corresponding asymptotic results for our generalisation of the FPOS distribution.

While the moments of the generalised distribution are complicated in general, there are some useful things that can be said about the asymptotic properties of the order statistics. Suppose we are willing to stipulate that the limiting empirical distribution of $z_1, z_2, z_3, \ldots$ is a continuous distribution function $F_Z$ with the corresponding quantile function $Q_Z$. Using the asymptotic results in Theorems 5-7 we can easily establish convergence in probability $\widetilde{U}_{(k)} \to \phi$ which then gives convergence in probability $Z_{(k)} \to Q_Z(\phi)$. This means that the order statistic generated from SRSWOR in a finite population will converge to the true quantile of the distribution under the limit condition we have previously specified. This is an intuitively reasonable property, given that SRSWOR from $z_1, \ldots, z_N$ becomes indistinguishable from IID sampling from the distribution $F_Z$ when we take the limit $N \to \infty$ (with other conditions on the limit as previously specified).

## 7. Conclusion

We have conducted a detailed examination of the distribution of order statistics arising in the case of simple-random-sampling-without-replacement (SRSWOR) from a finite population of values. This distribution is an interesting variation of the beta-binomial distribution using a location shift and re-parameterisation (we have called it the finite-population-order-statistic distribution using the acronym FPOS for short). We have examined some useful properties of the distribution including its moments, mixture characterisation and asymptotic properties.

Examination of the beta-binomial mixture representation shows a close connection between the scaled order statistics (scaled to be in the interior of the unit interval) and the distribution of the order statistics arising from IID sampling from the continuous uniform distribution on the unit interval. This provides an interesting connection between the distribution of order statistics is SRSWOR and the distribution of order statistics in continuous uniform sampling. Our analysis of the joint distribution of order statistics likewise shows a close connection, this time using the Dirichlet-multinomial mixture representation. One useful aspect of this result is that it allows direct pseudo-random generation of order statistics from SRSWOR without having to undertake the sampling and sorting inherent in standard generative methods.

Our analysis has also examined the properties of the distribution with respect to inferences for the population size $N$. We have established that the highest observed order statistic is the complete sufficient statistic for $N$ and all lower order statistics are conditionally ancillary for $N$. This confirms the intuition that inference about the population size ought to be based on only the highest observed order statistic, with any lower order statistics contributing no more information. This finding gives rise to a generalisation of the German tank problem, where we estimate the population size based on observation of an arbitrary order statistic. We have derived the properties of a generalised classical estimator of the population size and we have given a heuristic explanation of this estimator based on the idea of considering the observed order statistics as being equally spaced in the interior of the interval $[0, N+1]$. We have also shown how one can incorporate prior information to obtain the posterior distribution for the population size in a Bayesian analysis.

We hope that the present paper sheds some light on the distribution of order statistics under SRSWOR from a finite population. This is a common method of sampling and it is useful to see the distribution of the order statistics in this case. It is particularly interesting to see the connection to the distribution of order statistics in other well-known cases.

# Appendix: Proofs of Theorems and Computation Result

**PROOF OF THEOREM 1:** We have:

$$\mathbb{E}(X_{(k)}^{(r)}) = \sum_{x=k}^{N-n+k} x^{(r)} \cdot \text{FPOS}(x|k,n,N)$$

$$= \sum_{x=k}^{N-n+k} x^{(r)} \cdot \frac{\binom{x-1}{k-1}\binom{N-x}{n-k}}{\binom{N}{n}}$$

$$= \sum_{x=k}^{N-n+k} k^{(r)} \cdot \frac{\binom{x+r-1}{k+r-1}\binom{N-x}{n-k}}{\binom{N}{n}}$$

$$= \sum_{x=k}^{N-n+k} \frac{(N+1)^{(r)} \cdot k^{(r)}}{(n+1)^{(r)}} \cdot \frac{\binom{x+r-1}{k+r-1}\binom{N-x}{n-k}}{\binom{N+r}{n+r}}$$

$$= \frac{(N+1)^{(r)} \cdot k^{(r)}}{(n+1)^{(r)}} \sum_{x=k}^{N-n+k} \text{FPOS}(x+r|k+r, n+r, N+r)$$

$$= \frac{(N+1)^{(r)} \cdot k^{(r)}}{(n+1)^{(r)}},$$

which was to be shown. ∎

**PROOF OF THEOREM 2:** Substituting $r=1$ and $r=2$ in Theorem 1 gives the moments:

$$\mathbb{E}(X_{(k)}) = \frac{N+1}{n+1} \cdot k \qquad \mathbb{E}(X_{(k)}(X_{(k)}+1)) = \frac{(N+1)(N+2)}{(n+1)(n+2)} \cdot k(k+1).$$

The first result gives us the mean, and to obtain the variance we have:

$$\mathbb{V}(X_{(k)}) = \mathbb{E}(X_{(k)}^2) - \mathbb{E}(X_{(k)})^2$$

$$= \mathbb{E}(X_{(k)}(X_{(k)}+1)) - \mathbb{E}(X_{(k)}) - \mathbb{E}(X_{(k)})^2$$

$$= \frac{(N+1)(N+2)}{(n+1)(n+2)} \cdot k(k+1) - \frac{N+1}{n+1} \cdot k - \frac{N+1}{n+1} \cdot \frac{N+1}{n+1} \cdot k^2$$

$$= \frac{(N+1)k}{(n+1)^2(n+2)} \cdot \begin{bmatrix} (N+2)(n+1)(k+1) \\ -(n+1)(n+2) \\ -(N+1)(n+2)k \end{bmatrix}$$

$$= \frac{(N+1)k}{(n+1)^2(n+2)} \cdot \begin{bmatrix} (Nnk+Nn+Nk+N) \\ +(2nk+2n+2k+2) \\ -(n^2+3n+2) \\ -(Nnk+2Nk+nk+2k) \end{bmatrix}$$

$$= \frac{(N+1)k}{(n+1)^2(n+2)} \cdot [Nn - Nk + N + nk - n^2 - n]$$

$$= \frac{(N+1)(N-n)}{(n+1)^2(n+2)} \cdot k(n-k+1),$$

which was to be shown. ∎

**PROOF OF THEOREM 3:** The proof of the skewness and kurtosis formulae follows the same general method as in Theorem 2 —i.e., we substitute $r = 3$ and $r = 4$ into the factorial moment formula in Theorem 1 and then perform the necessary algebra to get the third and fourth central moments (and from these the skewness and kurtosis). The working is extremely cumbersome and so it is omitted here. In order to confirm correctness of the moment formulae, we compared the results calculated from these formulae to the moments produced by summing powers of deviations from the mean, multiplies by the mass values for the FPOS distribution. ∎

**PROOF OF THEOREM 4:** By application of the beta integral we have:

$$H(x|k, n, N) \equiv \int_0^1 \text{Bin}(x - k|N - n, u) \cdot \text{Beta}(u|k, n - k + 1)\, du$$

$$= \int_0^1 \binom{N-n}{x-k} \cdot u^{x-k}(1-u)^{N-n-x+k} \cdot \frac{n!}{(k-1)!\,(n-k)!} \cdot u^{k-1}(1-u)^{n-k} du$$

$$= \binom{N-n}{x-k} \cdot \frac{n!}{(k-1)!\,(n-k)!} \int_0^1 u^{x-1}(1-u)^{N-x} du$$

$$= \frac{(N-n)!}{(x-k)!\,(N-n-x+k)!} \cdot \frac{n!}{(k-1)!\,(n-k)!} \cdot \frac{(x-1)!\,(N-x)!}{N!}$$

$$= \frac{(N-n)!\,n!}{N!} \cdot \frac{(x-1)!}{(k-1)!\,(x-k)!} \cdot \frac{(N-x)!}{(n-k)!\,(N-n-x+k)!}$$

$$= \frac{\binom{x-1}{k-1}\binom{N-x}{n-k}}{\binom{N}{n}},$$

which was to be shown. ∎

**PROOF OF THEOREM 5:** The asymptotic forms and their asymptotic equivalence to the actual moments is easily shown by substitution of $n = \lambda N$ and $k = \phi \lambda N$ and cancelling of relevant terms for large $N$. Taking $N \to \infty$ then shows that the distribution is asymptotically unskewed and mesokurtic. ∎

**LEMMA 1:** Define the random function $H$ and the deterministic function $H_*$ by:

$$H(t) \equiv \frac{(N-n)U_{(k)} - (\mathbb{E}(X_{(k)}) - k)}{\mathbb{S}(X_{(k)})} \cdot t + \frac{(N-n)U_{(k)}(1 - U_{(k)})}{\mathbb{V}(X_{(k)})} \cdot \frac{t^2}{2},$$

$$H_*(t) \equiv \frac{\mathbb{C}\mathrm{ov}(U_{(k)}, U_{(k)}(1 - U_{(k)}))}{\mathbb{S}(X_{(k)})^3/(N-n)^2} \cdot \frac{t^3}{2} + \frac{\mathbb{V}(U_{(k)}(1 - U_{(k)}))}{\mathbb{S}(X_{(k)})^4/(N-n)^2} \cdot \frac{t^4}{4}.$$

At any point $t \in \mathbb{R}$ we have:

$$\mathbb{E}(H(t)) = \frac{n+1}{N+1} \cdot \frac{t^2}{2} \qquad \mathbb{V}(H(t)) = \frac{N-n}{N+1} \cdot t^2 + H_*(t).$$

Under the limit stipulated in Theorem 6 we have:

$$\lim_{N,n,k} H_*(t) = 0 \qquad \text{for all } t \in \mathbb{R}.$$

**PROOF OF LEMMA 1:** Using the raw moments of the beta distribution we have:

$$\mathbb{E}(U_{(k)}^r) = \frac{k}{n+1} \cdot \frac{k+1}{n+2} \cdot \ldots \cdot \frac{k+r-1}{n+r}.$$

With a bit of algebra, we can use this result to establish that:

$$\mathbb{E}(U_{(k)}) = \frac{k}{n+1},$$

$$\mathbb{V}(U_{(k)}) = \frac{k(n-k+1)}{(n+1)^2(n+2)},$$

$$\mathbb{E}(U_{(k)}(1 - U_{(k)})) = \frac{k(n-k+1)}{(n+1)(n+2)},$$

$$\mathbb{V}(U_{(k)}(1 - U_{(k)})) = \frac{k\begin{pmatrix} -n^4 - 10n^3 - 29n^2 \\ +kn^3 + 19kn^2 \\ -4k^3n - 10k^3 + 4k^2n^2 \\ +4k^2n - 12k^2 \\ +52kn + 34k - 32n - 12 \end{pmatrix}}{(n+1)^2(n+2)^2(n+3)(n+4)},$$

$$\mathbb{C}\mathrm{ov}(U_{(k)}, U_{(k)}(1 - U_{(k)})) = \frac{k(n-k+1)(n-2k+1)}{(n+1)^2(n+2)(n+3)}.$$

Using these moments and the relevant moments of $X_{(k)}$ in Theorem 2 we then have:

$$\mathbb{E}\left(\frac{(N-n)U_{(k)} - (\mathbb{E}(X_{(k)}) - k)}{\mathbb{S}(X_{(k)})}\right) = \frac{\mathbb{E}((N-n)U_{(k)} - (\mathbb{E}(X_{(k)}) - k))}{\mathbb{S}(X_{(k)})}$$

$$= \frac{1}{\mathbb{S}(X_{(k)})}\left[(N-n) \cdot \frac{k}{n+1} - \frac{N+1}{n+1} \cdot k + k\right]$$

$$= \frac{1}{\mathbb{S}(X_{(k)})} \cdot \frac{k}{n+1} [(N-n) - (N+1) + (n+1)]$$

$$= \frac{1}{\mathbb{S}(X_{(k)})} \cdot \frac{k}{n+1} \times 0$$

$$= 0,$$

$$\mathbb{E}\left(\frac{(N-n)U_{(k)}(1-U_{(k)})}{\mathbb{V}(X_{(k)})}\right) = \frac{N-n}{\mathbb{V}(X_{(k)})} \cdot \mathbb{E}\left(U_{(k)}(1-U_{(k)})\right)$$

$$= \frac{N-n}{\mathbb{V}(X_{(k)})} \cdot \frac{k(n-k+1)}{(n+1)(n+2)}$$

$$= \frac{n+1}{N+1},$$

$$\mathbb{V}\left(\frac{(N-n)U_{(k)}}{\mathbb{S}(X_{(k)})}\right) = \frac{(N-n)^2}{\mathbb{V}(X_{(k)})} \cdot \mathbb{V}(U_{(k)})$$

$$= \frac{(N-n)^2}{\mathbb{V}(X_{(k)})} \cdot \frac{k(n-k+1)}{(n+1)^2(n+2)}$$

$$= \frac{(N-n)^2}{\mathbb{V}(X_{(k)})} \cdot \frac{k(n-k+1)}{(n+1)^2(n+2)}$$

$$= \frac{N-n}{N+1}.$$

We therefore have:

$$\mathbb{E}(H(t)) = \mathbb{E}\left(\frac{(N-n)U_{(k)} - (\mathbb{E}(X_{(k)}) - k)}{\mathbb{S}(X_{(k)})}\right) \cdot t + \mathbb{E}\left(\frac{(N-n)U_{(k)}(1-U_{(k)})}{\mathbb{V}(X_{(k)})}\right) \cdot \frac{t^2}{2}$$

$$= \frac{n+1}{N+1} \cdot \frac{t^2}{2},$$

and:

$$\mathbb{V}(H(t)) = \mathbb{V}\left(\frac{(N-n)U_{(k)}}{\mathbb{S}(X_{(k)})}\right) \cdot t^2 + \mathbb{V}\left(\frac{(N-n)U_{(k)}(1-U_{(k)})}{\mathbb{V}(X_{(k)})}\right) \cdot \frac{t^4}{4}$$

$$+ \mathbb{C}\mathrm{ov}\left(\frac{(N-n)U_{(k)}}{\mathbb{S}(X_{(k)})}, \frac{(N-n)U_{(k)}(1-U_{(k)})}{\mathbb{V}(X_{(k)})}\right) \cdot \frac{t^3}{2}$$

$$= \frac{N-n}{N+1} \cdot t^2 + \frac{\mathbb{C}\mathrm{ov}(U_{(k)}, U_{(k)}(1-U_{(k)}))}{\mathbb{S}(X_{(k)})^3/(N-n)^2} \cdot \frac{t^3}{2} + \frac{\mathbb{V}(U_{(k)}(1-U_{(k)}))}{\mathbb{S}(X_{(k)})^4/(N-n)^2} \cdot \frac{t^4}{4}$$

$$= \frac{N-n}{N+1} \cdot t^2 + H_*(t).$$

This establishes the first part of the lemma, so it remains only to establish the limiting result. To do this, we first note that the limit stipulated in Theorem 6 implies:[1]

$$\text{Term}_1 = \frac{\mathbb{Cov}(U_{(k)}, U_{(k)}(1-U_{(k)}))}{\mathbb{S}(X_{(k)})^3/(N-n)^2} \cdot \frac{t^3}{2}$$

$$= (N-n)^2 \cdot \frac{k(n-k+1)(n-2k+1)}{(n+1)^2(n+2)(n+3)} \left(\frac{(n+1)^2(n+2)}{(N+1)(N-n)k(n-k+1)}\right)^{3/2} \cdot \frac{t^3}{2}$$

$$= \frac{(N-n)^{1/2}(n+1)^{1/2}}{(N+1)^{3/2}} \cdot \frac{n-2k+1}{n+3} \cdot \left(\frac{(n+1)(n+2)}{k(n-k+1)}\right)^{1/2} \cdot \frac{t^3}{2}$$

$$\simeq \frac{1-2\lambda}{N^{1/2}} \cdot \sqrt{\frac{\phi(1-\phi)}{\lambda(1-\lambda)}} \cdot \frac{t^3}{2} = \mathcal{O}(N^{-1/2}) \to 0.$$

$$\text{Term}_2 = \frac{\mathbb{V}(U_{(k)}(1-U_{(k)}))}{\mathbb{S}(X_{(k)})^4/(N-n)^2} \cdot \frac{t^4}{4}$$

$$= \frac{k\begin{pmatrix} -n^4 - 10n^3 - 29n^2 \\ +kn^3 + 19kn^2 \\ -4k^3n - 10k^3 + 4k^2n^2 \\ +4k^2n - 12k^2 \\ +52kn + 34k - 32n - 12 \end{pmatrix}}{(n+1)^2(n+2)^2(n+3)(n+4)} \cdot \left(\frac{(n+1)^2(n+2)}{(N+1)(N-n)k(n-k+1)}\right)^2 \cdot \frac{t^4}{4}$$

$$= \frac{\begin{pmatrix} -n^4 - 10n^3 - 29n^2 \\ +kn^3 + 19kn^2 \\ -4k^3n - 10k^3 + 4k^2n^2 \\ +4k^2n - 12k^2 \\ +52kn + 34k - 32n - 12 \end{pmatrix}}{(N+1)^2(N-n)^2} \cdot \frac{(n+1)^2}{(n+3)(n+4)k(n-k+1)^2} \cdot \frac{t^4}{4}$$

$$\simeq -\frac{\phi^4 + \lambda\phi^3 - 4\lambda^3\phi + 4\lambda^2\phi^2}{\phi^2(1-\phi)^2} \cdot \frac{1}{\lambda(1-\lambda)^2} \cdot \frac{1}{N^2} \cdot \frac{t^4}{4} = \mathcal{O}(N^{-2}) \to 0.$$

From these results we therefore have $H_*(t) = \text{Term}_1 + \text{Term}_2 \to 0$ which establishes the limit and completes the proof. ∎

**LEMMA 2:** Define the random function $G$ by:

$$G(t, \xi) \equiv \frac{(N-n)\xi U_{(k)}(1-\xi U_{(k)})(½ - \xi U_{(k)})}{\mathbb{S}(X_{(k)})^3} \cdot \frac{t^3}{3} \qquad \text{for } t \in \mathbb{R}, 0 \leq \xi \leq 1.$$

Under the limit stipulated in Theorem 6 we have:

---

[1] The operator $\simeq$ refers to asymptotic equivalence; in much literature this is denoted as $\sim$ but we have avoided that latter notation since it is already used in this paper to refer to distributional equivalence. Formally, we say that $f(N,n,k) \simeq g(N,n,k)$ under the stipulated limit if $\lim_{N,n,k} f(N,n,k)/g(N,n,k) = 1$.

$$\plim_{N,n,k} G(t,\xi) = 0 \qquad \text{for all } t \in \mathbb{R}, 0 \leq \xi \leq 1.$$

**PROOF OF LEMMA 2:** Since $0 \leq \xi U_{(k)} \leq 1$ we have:

$$|\xi U_{(k)}(1 - \xi U_{(k)})(\tfrac{1}{2} - \xi U_{(k)})| = |\xi U_{(k)}||1 - \xi U_{(k)}||\tfrac{1}{2} - \xi U_{(k)}| \leq 1.$$

Using the relevant moments of $X_{(k)}$ in Theorem 2 we then have:

$$\plim_{N,n,k}|G(t,\xi)| = \plim_{N,n,k} \left| \frac{(N-n)\xi U_{(k)}(1-\xi U_{(k)})(\tfrac{1}{2} - \xi U_{(k)})}{\mathbb{S}(X_{(k)})^3} \cdot \frac{t^3}{3} \right|$$

$$\leq \lim_{N,n,k} \left| \frac{(N-n)}{\mathbb{S}(X_{(k)})^3} \cdot \frac{t^3}{3} \right|$$

$$= \lim_{N,n,k} \left| \frac{(n+1)^3(n+2)^{3/2}}{(N+1)^{3/2}(N-n)^{1/2}k^{3/2}(n-k+1)^{3/2}} \cdot \frac{t^3}{3} \right|$$

$$= \lim_{N,n,k} \left| \frac{\phi^{3/2}}{(1-\phi)^{1/2}\lambda^{3/2}(1-\lambda)^{3/2}} \cdot \frac{t^3}{3} \cdot \frac{1}{\sqrt{N}} \right| = 0.$$

$$= \frac{\phi^{3/2}}{(1-\phi)^{1/2}\lambda^{3/2}(1-\lambda)^{3/2}} \cdot \frac{t^3}{3} \cdot \lim_{N,n,k} \frac{1}{\sqrt{N}} = 0.$$

(Note that the penultimate step uses the asymptotic equivalence in the stipulated limit.) Since $\plim_{N,n,k}|G(t,\xi)| = 0$ we have $\plim_{N,n,k} G(t,\xi) = 0$ which was to be shown. ∎

**PROOF OF THEOREM 6:** To facilitate our analysis we will use the function $R: (0,1) \times \mathbb{R} \to \mathbb{R}$ defined by:

$$R(a,t) \equiv \log(1 - a(1 - e^t)).$$

The second-order Maclaurin representation for this function (with remainder) is:

$$R(a,t) = a \cdot t + a(1-a) \cdot \frac{t^2}{2} + \xi a(1-\xi a)(\tfrac{1}{2} - \xi a) \cdot \frac{t^3}{3},$$

where $0 \leq \xi \leq 1$ is the proportion used in the remainder term. Now, recall the representation:

$$X_{(k)} \sim k + \text{Bin}(N-n, U_{(k)}) \qquad U_{(k)} \sim \text{Beta}(k, n-k+1).$$

Combining this result with the known form of the moment generating function for the binomial distribution, and using the functions $H$ and $G$ in Lemmas 1-2, we get:

$$m_*(t, U_{(k)}) \equiv \mathbb{E}\left(\exp\left(t \cdot \frac{X_{(k)} - \mathbb{E}(X_{(k)})}{\mathbb{S}(X_{(k)})}\right) \Big| U_{(k)}\right)$$

$$= \mathbb{E}\left(\exp\left(t \cdot \frac{(X_{(k)} - k) - (\mathbb{E}(X_{(k)}) - k)}{\mathbb{S}(X_{(k)})}\right) \Big| U_{(k)}\right)$$

$$= \exp\left(-\frac{\mathbb{E}(X_{(k)}) - k}{\mathbb{S}(X_{(k)})} \cdot t\right) \mathbb{E}\left(\exp\left(\frac{t}{\mathbb{S}(X_{(k)})} \cdot (X_{(k)} - k)\right) \bigg| U_{(k)}\right)$$

$$= \exp\left(-\frac{\mathbb{E}(X_{(k)}) - k}{\mathbb{S}(X_{(k)})} \cdot t\right)\left(1 - U_{(k)}\left(1 - \exp\left(\frac{t}{\mathbb{S}(X_{(k)})}\right)\right)\right)^{N-n}$$

$$= \exp\left(-\frac{\mathbb{E}(X_{(k)}) - k}{\mathbb{S}(X_{(k)})} \cdot t + (N-n) \cdot R\left(U_{(k)}, \frac{t}{\mathbb{S}(X_{(k)})}\right)\right)$$

$$= \exp(H(t) + G(t, \xi)).$$

In Lemma 2 we established that $G(t, \xi) \to 0$ in probability, so we have:

$$\plim_{N,n,k} m_*(t, U_{(k)}) = \plim_{N,n,k} \exp(H(t)).$$

Now, for any $t \in \mathbb{R}$ the random variable $H(t)$ is a quadratic form of a beta random variable, which converges to a quadratic form of a normal random variable, which in turn converges in distribution to the normal distribution. This means that —under the stipulated limit— the random variable $\exp(H(t))$ converges in distribution to the lognormal distribution. Using the expected value of the lognormal distribution we then have:

$$\lim_{N,n,k} m(t) = \lim_{N,n,k} \mathbb{E}(m_*(t, U_{(k)}))$$

$$= \lim_{N,n,k} \mathbb{E}\left(\plim_{N,n,k} m_*(t, U_{(k)})\right)$$

$$= \lim_{N,n,k} \mathbb{E}\left(\plim_{N,n,k} \exp(H(t))\right)$$

$$= \lim_{N,n,k} \exp\left(\mathbb{E}(H(t)) + \frac{1}{2} \cdot \mathbb{V}(H(t))\right)$$

$$= \lim_{N,n,k} \exp\left(\frac{n+1}{N+1} \cdot \frac{t^2}{2} + \frac{N-n}{N+1} \cdot \frac{t^2}{2} + \frac{1}{2} \cdot H_*(t)\right)$$

$$= \lim_{N,n,k} \exp\left(\frac{t^2}{2} + \frac{1}{2} \cdot H_*(t)\right)$$

$$= \exp\left(\frac{t^2}{2} + \frac{1}{2} \cdot \lim_{N,n,k} H_*(t)\right)$$

$$= \exp\left(\frac{t^2}{2}\right),$$

where the fourth-last and second-last steps both follow from results in Lemma 1. This establishes the result to be shown. ∎

**PROOF OF THEOREM 7:** See Khan (1994). ∎

**PROOF OF THEOREM 8:** By application of the multivariate beta integral we have:

$$H(\boldsymbol{x}_*|\boldsymbol{k}_*, n, N) \equiv \int_\Theta \text{Mu}(\Delta_x - \Delta_k | N - n, \boldsymbol{s}) \cdot \text{Dirichlet}(\boldsymbol{s}|\Delta_k) \, d\boldsymbol{s}$$

$$= \int_\Theta \Gamma(N - n + 1) \prod_{i=1}^{r+1} \frac{s_i^{\Delta_{x,i} - \Delta_{k,i}}}{\Gamma(\Delta_{x,i} - \Delta_{k,i} + 1)} \cdot \Gamma(n + 1) \prod_{i=1}^{r+1} \frac{s_i^{\Delta_{k,i} - 1}}{\Gamma(\Delta_{k,i})} \, d\boldsymbol{s}$$

$$= \Gamma(N - n + 1)\Gamma(n + 1) \int_\Theta \prod_{i=1}^{r+1} \frac{s_i^{\Delta_{x,i} - 1}}{\Gamma(\Delta_{x,i} - \Delta_{k,i} + 1)\Gamma(\Delta_{k,i})} \, d\boldsymbol{s}$$

$$= \frac{\Gamma(N - n + 1)\Gamma(n + 1)}{\prod_i \Gamma(\Delta_{x,i} - \Delta_{k,i} + 1)\Gamma(\Delta_{k,i})} \int_\Theta \prod_{i=1}^{r+1} s_i^{\Delta_{x,i} - 1} \, d\boldsymbol{s}$$

$$= \frac{\Gamma(N - n + 1)\Gamma(n + 1)}{\prod_i \Gamma(\Delta_{x,i} - \Delta_{k,i} + 1)\Gamma(\Delta_{k,i})} \cdot \frac{\prod_i \Gamma(\Delta_{x,i})}{\Gamma(N + 1)}$$

$$= \frac{\Gamma(N - n + 1)\Gamma(n + 1)}{\Gamma(N + 1)} \prod_{i=1}^{r+1} \frac{\Gamma(\Delta_{x,i})}{\Gamma(\Delta_{x,i} - \Delta_{k,i} + 1)\Gamma(\Delta_{k,i})}$$

$$= \frac{1}{\binom{N}{n}} \cdot \prod_{i=1}^{r+1} \binom{\Delta_{x,i} - 1}{\Delta_{k,i} - 1}$$

$$= \frac{\binom{x_1 - 1}{k_1 - 1}\binom{x_2 - x_1 - 1}{k_2 - k_1 - 1}\binom{x_3 - x_2 - 1}{k_3 - k_2 - 1} \cdots \binom{x_r - x_{r-1} - 1}{k_r - k_{r-1} - 1}\binom{N - x_r}{n - k_r}}{\binom{N}{n}},$$

which was to be shown. ∎

**LEMMA 3:** Given any function $g: \mathbb{N} \to \mathbb{R}$ we define the corresponding function:

$$G(N) \equiv \sum_{x=k}^{N-n+k} g(x) \cdot \frac{(N - x)! \, (x - 1)!}{(N - n - x + k)! \, (x - k)!}.$$

If $G(N) = 0$ for all $N \geq n$ then $g(x) = 0$ for all $x \geq k$.

**PROOF OF LEMMA 3:** We prove this lemma using strong induction on $N$. For all remaining steps we assume that the antecedent condition of the lemma is true. Taking $N = n$ gives:

$$0 = G(n) = \sum_{x=k}^{k} g(x) \cdot \frac{(n - x)! \, (x - 1)!}{(k - x)! \, (x - k)!}$$

$$= g(k) \cdot (n - k)! \, (k - 1)!$$

which implies that $g(k) = 0$. This gives the base case for the induction. Now, suppose that $g(k) = \cdots = g(k + r - 1) = 0$ for some $r \geq 1$. Taking $N = n + r$ gives:

$$0 = G(n+r) = \sum_{x=k}^{r+k} g(x) \cdot \frac{(n+r-x)!\,(x-1)!}{(r-x+k)!\,(x-k)!}$$

$$= g(r+k) \cdot \frac{(n-k)!\,(r+k-1)!}{r!}$$

which implies that $g(k + r) = 0$. This establishes the induction step, which is sufficient to establish that $g(x) = 0$ for all $x \geq k$. ∎

**PROOF OF THEOREM 9(a):** Taking $\boldsymbol{k}_{**} = (k_1, \ldots, k_{r-1})$ and $\boldsymbol{x}_{**} = (x_1, \ldots, x_{r-1})$ we write:

$$\text{FPOS}(\boldsymbol{x}_*|\boldsymbol{k}_*, n, N) = \text{FPOS}(\boldsymbol{x}_{**}|\boldsymbol{k}_{**}, k_r - 1, x_r - 1) \cdot \text{FPOS}(x_r|k_r, n, N).$$

Since only the latter part depends on $N$ this gives the Fisher-Neyman factorisation of the joint distribution, which shows that $X_{(k_r)}$ is a sufficient statistic for $N$. To prove completeness, we take an arbitrary function $g: \mathbb{N} \to \mathbb{R}$ and set $G: \mathbb{N} \to \mathbb{R}$ to be the corresponding function defined in Lemma 3. We then have:

$$\mathbb{E}(g(X_{(k_r)})) = \sum_x g(x) \cdot \text{FPOS}(x|k_r, n, N)$$

$$= \sum_{x=k}^{N-n+k} g(x) \cdot \frac{\binom{x-1}{k-1}\binom{N-x}{n-k}}{\binom{N}{n}}$$

$$= \frac{(N-n)!\,n!}{N!\,(n-k)!\,(k-1)!} \sum_{x=k}^{N-n+k} g(x) \cdot \frac{(N-x)!\,(x-1)!}{(N-n-x+k)!\,(x-k)!}$$

$$= \frac{(N-n)!\,n!}{N!\,(n-k)!\,(k-1)!} \cdot G(N).$$

Since $1 \leq k \leq n \leq N$ the multiplicative term at the front of this expression is strictly positive. Consequently, we have $\mathbb{E}(g(X_{(k_r)})) = 0$ if and only if $G(N) = 0$. Applying Lemma 3 ensures that if $\mathbb{E}(g(X_{(k_r)})) = 0$ for all $N \geq n$ then we have $g(x) = 0$ for all $x \geq k$, which establishes that $X_{(k_r)}$ is complete with respect to $N$. This establishes that $X_{(k_r)}$ is a complete sufficient statistic for $N$ which was to be shown. ∎

**PROOF OF THEOREM 9(b):** We again take $\boldsymbol{k}_{**} = (k_1, \ldots, k_{r-1})$ and $\boldsymbol{x}_{**} = (x_1, \ldots, x_{r-1})$. Using this notation the distribution of $\boldsymbol{X}_{**}$ conditional on $X_{(k_r)}$ has probability mass function:

$$\mathbb{P}(\boldsymbol{X}_{**} = \boldsymbol{x}_{**}|X_{(k_r)} = x_r) = \frac{\text{FPOS}(\boldsymbol{x}_*|\boldsymbol{k}_*, n, N)}{\text{FPOS}(x_r|k_r, n, N)}$$

$$= \frac{\binom{x_1-1}{k_1-1}\binom{x_2-x_1-1}{k_2-k_1-1}\binom{x_3-x_2-1}{k_3-k_2-1}\cdots\binom{x_r-x_{r-1}-1}{k_r-k_{r-1}-1}}{\binom{x_r-1}{k_r-1}}.$$

Since this distribution does not depend on $N$ this means that $\boldsymbol{X}_{**}$ is ancillary for $N$ conditional on $X_{(k_r)}$, which was to be shown. ∎